# WAITING FOR REGULATORY SEQUENCES TO APPEAR

By Richard Durrett[1] and Deena Schmidt[2]

*Cornell University*

One possible explanation for the substantial organismal differences between humans and chimpanzees is that there have been changes in gene regulation. Given what is known about transcription factor binding sites, this motivates the following probability question: given a 1000 nucleotide region in our genome, how long does it take for a specified six to nine letter word to appear in that region in some individual? Stone and Wray [*Mol. Biol. Evol.* **18** (2001) 1764–1770] computed 5,950 years as the answer for six letter words. Here, we will show that for words of length 6, the average waiting time is 100,000 years, while for words of length 8, the waiting time has mean 375,000 years when there is a 7 out of 8 letter match in the population consensus sequence (an event of probability roughly 5/16) and has mean 650 million years when there is not. Fortunately, in biological reality, the match to the target word does not have to be perfect for binding to occur. If we model this by saying that a 7 out of 8 letter match is good enough, the mean reduces to about 60,000 years.

**1. Introduction.** At a genetic level humans and chimpanzees are closely related, with 98.7% of their DNA identical. It has long been speculated that many of the obvious differences between the two species are due to changes in regulatory sequences that control how genes are expressed (King and Wilson [15]). A regulatory sequence is a short sequence of DNA (in vertebrates many are 6–9 nucleotides long) which is a binding site for transcription factors that promote or inhibit transcription of the DNA to make proteins. It does not take an advanced knowledge of biology to appreciate the importance of gene regulation. All of the cells in our body have the same 3 billion nucleotide instruction set, but each of the tissues in our body requires its own specialized set of proteins.

Received December 2005; revised July 2006.

[1]Supported in part by an NSF/NIGMS grant.

[2]Supported in part by an NSF graduate fellowship.

Comparison of experimentally verified promoters in 51 genes by Dermitzakis and Clark [7], showed that roughly 1/3 of sites functional in humans are not functional in rodents. This shows that in the 90 million years since the divergence of humans and rodents, there has been a significant evolution of transcription factor binding sites. However, there are only a few documented cases of human specific regulatory evolution since our divergence from the other primate; see [14] and [19] for recent examples and references to earlier work.

These observations suggest the question: is the evolution of regulatory sequences sufficiently rapid to contribute to the differences between humans and chimpanzees? To begin to turn this into a probability problem, we note that regulatory sequences occur in the upstream region of a gene (i.e., before the gene if we are reading the DNA strand in the order in which it is transcribed). They are typically within 1 kb (kilobase = 1000 nucleotides) of the beginning of the gene. Our probability question then is: given a word of length $W$ from the DNA alphabet $\{A, C, G, T\}$ and a 1 kb region, how long do we have to wait for mutation to create this word in a given 1 kb region in the genome of some human? In formulating the question, we are assuming that once created, the new mutation will confer a substantial benefit and, hence, will with high probability spread rapidly through the population. We will return to this point in Section 2.5.

Stone and Wray [20] studied this problem by simulation. The numbers they give in the first column of page 1767 are not consistent with the values in their Table 1 for a 2 kb region. Working backward and rounding to simplify the arithmetic, we infer that they found that, for words of length $W = 6$, the waiting time for the appearance of a word in a 2 kb region of one individual required an average of 952 mutations or 0.476 mutations per nucleotide. They use $\mu = 10^{-9}$ as an estimate of the mutation rate per nucleotide per generation, so this translates into 476 million generations. To convert this into an estimate for the human population, they assumed that all individuals evolve independently and divided by $10^6$ individuals times 2 DNA strands (since humans are diploid) to arrive at an estimate of 238 generations. Using 25 years for the generation time of humans, they concluded that the desired new sequence would appear in an average of 5,950 years.

As MacArthur and Brookfield [17] have already pointed out, there are two problems with this computation. The first is that humans are closely related so the evolution of the DNA sequences of different individuals is not independent. Indeed, two randomly chosen individuals differ in 1 of 1000 nucleotides of their DNA. The second problem is that $10^6$ is a substantial overestimate for the "effective size" of the human population. To explain the last term, we note that in a homogeneously mixing population of $N$ diploid individuals, the genealogy of two copies of a nucleotide coalesces in each generation with probability $1/2N$ and each nucleotide experiences mutations with probability $2\mu$. Thus, the probability of a mutation before coalescence, which will make the two individuals different, is

$$\frac{2\mu}{1/2N + 2\mu} = \frac{4N\mu}{1 + 4N\mu}.$$

Using $\mu = 2.5 \times 10^{-8}$ for the mutation rate to simplify the arithmetic, setting the fraction equal to $1/1000$, approximating the denominator by 1, and solving gives $N = 10^4$, a number that is commonly used for the effective size of the human population. Using $10^4$ instead of $10^6$ in Stone and Wray's computation increases their estimate to 595,000 years. However, the far more substantial problem is the lack of independence. This further increases the time and brings up the question: is the evolution of regulatory sequences rapid enough to make a substantial contribution in the 6 million years since the divergence of humans and chimpanzees?

We build up to our solution (case 4 below) by considering the waiting time to find a prespecified $W$ letter DNA word in the following steps:

1. $W$ nucleotides in one DNA sequence.
2. A segment of $L$ nucleotides in one DNA sequence.
3. $W$ nucleotides in a population of $N$ diploid individuals.
4. A segment of $L$ nucleotides in a population of $N$ diploid individuals.

**2. Results.** In this section we will describe our results that give approximate answers for the four problems and their implications for biology. Before entering into the details, we recall some of Aldous' [1] thoughts on the philosophy of approximations, heuristics and limit theorems:

> "The proper business of probabilists is computing probabilities. Often exact calculations are tedious or impossible, so we resort to approximations. A limit theorem is an assertion of the form 'the error in a certain approximation tends to 0 as (say) $N \to \infty$.' Call such a limit theorem *naive* if there is no explicit bound in terms of $N$ and the parameters of the underlying process. Such theorems are so prevalent in theoretical and applied probability that people seldom stop to ask their purpose. Given a serious applied problem involving specific parameters, the natural first steps are to seek rough analytic approximations and to run computer simulations; the next step is to do careful numerical analysis. It is hard to give any argument for the relevance of a proof of a naive limit theorem, except as a vague reassurance that your approximation is sensible, and a good heuristic argument seems equally reassuring."

Throughout we will concentrate on DNA words of length $W = 6$ and $W = 8$, which are sizes appropriate for human transcription factor binding sites and have different qualitative behavior. For the first and simplest problem, we are able to clearly quantify the errors in our approximation. For the second problem, a result of Arratia, Goldstein and Gordon [3] bounds the error in an associated Poisson approximation. However, to translate this into a result about waiting times, and to study the more complex third and fourth problems, the best we have been able to do is to suggest approximations and explain why they should be accurate. Because of this, we label our results as approximations rather than theorems. The main problem is that since we are interested in small values of $W$, it is not sensible to investigate the limit $W \to \infty$. Perhaps others can succeed in finding more precise estimates of our errors. However, as we will see later, even our current level of understanding allows us to make important qualitative conclusions about the tempo of evolution of regulatory sequences. These have been mentioned in the abstract and will be explained in more detail in Section 2.5.

### 2.1. *W nucleotides in one DNA sequence.*

The simplest situation occurs when we consider waiting for a target $W$ letter word to appear at a specified $W$ nucleotides in one DNA sequence. To remove the waiting times between mutations, we consider a discrete time Markov chain $X_n$ that changes each time there is a mutation and gives the number of letters that match the $W$ letter target sequence. In biological reality not all mutations occur with equal probability, but for mathematical simplicity, we will suppose that they do. This simplification should not affect the order of magnitude of the estimate of our waiting time. $X_n$ has state space $S = \{0, 1, \ldots, W\}$ and transition probabilities

$$\begin{aligned}
p(x, x-1) &= x/W, \\
p(x, x+1) &= (1/3)(W - x)/W, \\
p(x, x) &= (2/3)(W - x)/W,
\end{aligned}$$

(1)

where all other $p(x, y) = 0$. To explain this, we note that the state decreases by one if we choose a matching letter to mutate. If we choose a nonmatching letter, the state increases by one if and only if we get the right one of the three possible mutations. In equilibrium, letters are random so the stationary distribution $\pi$ for $X_n$ is Binomial$(W, 1/4)$.

To understand the distribution of $T_W = \inf\{t : X_t = W\}$, we use Proposition 23 in Chapter 3 of [2] which implies that if we consider the continuous time chain $X_t$ that jumps at rate 1, and $A$ is any set, the hitting time $T_A = \inf\{t : X_t \in A\}$ satisfies

$$\sup_t |P_\pi(T_A > t) - \exp(-t/E_\pi T_A)| \leq \tau_2/E_\pi T_A,$$

(2)

where $\tau_2$ is the relaxation time, which Aldous and Fill define (see [2], page 19) to be 1 over the spectral gap. The continuous time chain $X_t$ has $\tau_2 = 3W/4$ (see the Appendix for details). As we will see in a moment, $E_\pi T_W \geq 4^W$, so the error is $<0.0011$ when $W = 6$ and $<0.0001$ when $W = 8$.

To compute the mean of $T_W$, we first consider starting the chain at the target sequence, that is, $P_W(X_0 = W) = 1$. Let $T_x^+ = \inf\{n \geq 1 : X_n = x\}$. A classic result of Kac (see, e.g., Theorem (3.3) in Chapter 6 of [9]) states:

TABLE 1

|                 | $W = 6$ | $W = 8$ |
|-----------------|---------|---------|
| $E_\pi T_W$     | 4,420   | 69,088  |
| $E_0 T_W$       | 4,431   | 69,104  |
| $4^W/(1-a)$     | 4,456   | 69,152  |

THEOREM 1. *If $X_n$ is an irreducible discrete time Markov chain on a finite state space with stationary distribution $\pi$, then $E_x T_x^+ = 1/\pi(x)$.*

In our case, $E_W T_W^+ = 1/\pi(W) = 4^W$. To compute $E_0 T_W$, we let $a = P_W(T_W^+ < T_0) = P_{W-1}(T_W < T_0)$. In the Appendix we show that

$$W \quad 6 \quad 8$$
$$a \quad 0.08093 \quad 0.05228.$$

Assuming that return times with $T_W^+ < T_0$ are small and dropping them from the expected value, we have

$$E_W T_W^+ \approx (1-a) E_0 T_W$$

or $E_0 T_W \approx 4^W/(1-a)$. Here and in what follows $\approx$ is read "is approximately" and has no precise mathematical meaning.

This is also the answer that comes from Aldous' Poisson clumping heuristic [1]. "Sparse random sets often resemble i.i.d. clumps centered at points of a Poisson process." Here the clumps consist of the returns to $W$ that come *soon* after the first one, which we make precise as returns to $W$ before hitting 0. Once the chain hits $W$, then it will visit $W$ for a geometric number of times with mean $1/(1-a)$ before it hits 0, so to get the correct mean, we should multiply $1/\pi(W)$ by the expected clump size $EC = 1/(1-a)$.

Since $P_\pi(T_W < T_0)$ is small, we expect that $E_\pi T_W \approx E_0 T_W \approx 4^W/(1-a)$. As Table 1 shows, the error in the approximation is less than 1% when $W = 6$ and less than 0.1% when $W = 8$.

The derivation of the numbers in the table can be found in the Appendix, where we compute other quantities that we need associated with the mutation chain.

APPROXIMATION 1. *For the continuous time mutation chain that jumps at rate 1, under $P_\pi$, $T_W$ is approximately exponential with mean $4^W/(1-a)$.*

To be precise, (2) shows that the distribution function $P_\pi(T_W > t)$ is uniformly within $3W/4^{W+1}$ of the exponential with mean $E_\pi T_W$, and our numerical computations show $E_\pi T_W$ differs from $4^W/(1-a)$ by less than 1%.

2.2. *A segment in one DNA sequence.* We next consider waiting for a target $W$ letter word to appear somewhere in a sequence of $L$ nucleotides in one DNA sequence. To be precise, we have a Markov chain $X_t \in \{A, C, G, T\}^L$ that jumps at rate $L$, and when it jumps a randomly chosen letter is changed to a different randomly chosen letter.

While $L$ is general, we are thinking about $L$ as being 1 or 2 thousand nucleotides. To avoid problems with the ends of the region, we will make the nonbiological assumption that our DNA region is a circle. It follows from symmetry that if we start from a random initial configuration, then when the word first appears, the probability it occurs in the circle but not in the flat DNA is $(W - 1)/L$, which is small for the values of $W$ and $L$ that we consider.

Let $A$ be the event that the target word is found somewhere in the sequence of $L$ nucleotides (with wrap around). This time (2) does not work well, since if each site jumps at rate 1, $\tau_2 = 3/4$ and as we will see, $E_\pi T_A \approx 4^W/(WL)$. Thus, the error bound is $(3WL)/4^{W+1}$. When $L = 1024 = 4^5$ this is

$$18/16 \text{ for } W = 6 \quad \text{and} \quad 6/64 \text{ for } W = 8.$$

As we will see in a minute, there is a good reason for a bad bound when $W = 6$. The waiting time starting from a random initial condition has an atom at 0 of size $\approx 1/4$, so the distribution cannot be well approximated by an exponential. It is for this reason we approach the waiting time problem using a version of the Chen–Stein method described in [3].

Let $I = \{0, 1, 2, \ldots, L - 1\}$. For $\alpha \in I$, let $A_\alpha = \{\alpha, \alpha + 1, \ldots, \alpha + W - 1\} \bmod L$ be the $W$ letter strip starting at index $\alpha$. Fix a time $T \geq 0$ and let $Y_\alpha$ be 1 if the target word appears in $A_\alpha$ by time $T$, and 0 otherwise. Note that the definition of $Y_\alpha$ depends on $T$ even though we have not recorded this dependence in the notation. Let $B_\alpha = \{\alpha - (W - 1), \ldots, \alpha + (W - 1)\} \bmod L$. If $\beta \in I \backslash B_\alpha$, then $Y_\alpha$ and $Y_\beta$ are independent, so we call $B_\alpha$ the neighborhood of dependence of $\alpha$.

Let $p_\alpha = P(Y_\alpha = 1)$, let $V = \sum_{\alpha \in I} Y_\alpha$ be the number of occurrences of the target word by time $T$ and $\lambda = EV$ be the expected value. Let $Z$ be a Poisson random variable with $EZ = \lambda$. We want to estimate the total variation distance $d_{\mathrm{TV}}(V, Z) = \sup_A |P(V \in A) - P(Z \in A)|$. Define

$$
\begin{aligned}
b_1 &= \sum_{\alpha \in I} \sum_{\beta \in B_\alpha} p_\alpha p_\beta, \\
b_2 &= \sum_{\alpha \in I} \sum_{\beta \in B_\alpha \backslash \{\alpha\}} E[Y_\alpha Y_\beta], \\
b_3 &= \sum_{\alpha \in I} E|E[Y_\alpha - p_\alpha | Y_\beta : \beta \notin B_\alpha]|.
\end{aligned}
\tag{3}
$$

In our case, $Y_\alpha$ is independent of $Y_\beta$, $\beta \notin B_\alpha$, so $b_3 = 0$ and Theorem 1 in [3] simplifies to

$$(4) \qquad d_{\text{TV}}(V, Z) \leq 2(b_1 + b_2)\left(\frac{1 - e^{-\lambda}}{\lambda}\right).$$

The last factor is always $\leq 1$, but can be very helpful for large $\lambda$.

a. *Success in the initial condition.* If we put down a random sequence of $L$ nucleotides, then the expected number of times the target word is present in a random initial configuration is $\gamma = L/4^W$. Taking $L = 1024 = 4^5$ for simplicity,

| $W$ | 6 | 8 |
|---|---|---|
| $\gamma = L/4^W$ | 0.25 | 0.015625 |
| $1 - e^{-\gamma}$ | 0.22119 | 0.015504. |

When $W = 8$, $L/4^W = 0.015625$ gives an upper bound on the probability of a match in the initial condition, so we will ignore this possibility in our approximation. We will now use the Arratia, Goldstein and Gordon [3] result stated in (4) to argue for the following:

APPROXIMATION 2a. *When $W = 6$ and $L = 1024$ if we exclude* 52 *repetitive words, then the number of matches in the initial condition is approximately Poisson with mean* 1/4, *with an error of at most* 0.0063 *in the total variation distance.*

When $T = 0$, $p_\alpha = 4^{-W}$ so

$$(5) \qquad b_1 = L(2W - 1)p_\alpha^2 = \gamma(2W - 1)/4^W,$$

which is $0.00268\gamma$ when $W = 6$. As we will see in Section 3, the bound on $b_2$ depends on the values of $k$ for which the word and its shift by $k$ letters agree exactly on the overlap. In Table 2 we have given results for the possible values of $k$ and an example pattern for each case.

Ignoring the first three categories, which account for 48 of the 4092 nonconstant words of length 6, the maximum value of $2(b_1 + b_2)$ is $0.02490\gamma = 0.006225$.

b. *Number of occurrences by time $T$.* Our next task is to use the Arratia, Goldstein and Gordon [3] result to find a Poisson approximation for the number of occurrences of the target word in a segment of length $L$ by time $T$. To evaluate the bound in (4), we note that, as in (5),

$$(6) \qquad b_1 = L(2W - 1)p_\alpha^2 = \lambda^2(2W - 1)/L,$$

TABLE 2

| $k$ | Pattern | $b_2/\gamma$ |
|------|---------|--------------|
| 2, 4 | ACACAC | 0.13281 |
| 3, 5 | ACAACA | 0.03320 |
| 3 | ACGACG | 0.03125 |
| 4, 5 | AACGAA | 0.00977 |
| 4 | ACGTAC | 0.00781 |
| 5 | ACGTCA | 0.00193 |
| none | ACGTAG | 0 |

TABLE 3

|        | tv       |          | tv       |
|--------|----------|----------|----------|
| AACCGT | 0.134229 | ACAGCTGT | 0.070616 |
| ACGCTA | 0.142948 | ACAAGGGC | 0.075011 |
| ACAGCA | 0.183293 | ACAGACAG | 0.100627 |
| AACGAA | 0.229230 | AAAAAACA | 0.145996 |
| ACAACA | 0.302622 | AACAACAA | 0.163626 |
| ACACAC | 0.465964 | ACACACAC | 0.337132 |

since $\lambda = Lp_\alpha$. Again the estimation of $b_2$ is more complicated than that of $b_1$ and depends on how the word overlaps with its shifts (see Section 3.2 for more details). Suppose, for concreteness, that $L = 1024$, $W = 6$ and $\lambda = 1$. In each case $b_1 = 0.010742$. We have computed the bounds on the total variation distance (tv) for all 4092 nonconstant words of length 6 and selected six words to illustrate the range of estimates of the total variation bound from the best $AACCGT$ to the worst $ACACAC$. We have also done the computation for the 65,532 nonconstant words of length 8. In each case $b_1 = 0.014648$. This time the word $ACAGCTGT$ achieves the best value, while $ACACACAC$ is again the worst (see Table 3).

Even our best result is not very comforting for someone who wants to use the result to estimate probabilities. To improve the quality of our approximation, we remember the Poisson clumping heuristic. We define new indicator variables $\bar{Y}_\alpha$ that only count a hit in strip $A_\alpha$ if there has been no hit in any overlapping strip $A_\beta$ for $\beta \in B_\alpha$ since the last time the number of matches in the strip $A_\alpha$ was 0. If two or more hits occur simultaneously, we pick one of the $\alpha$ at random to have $\bar{Y}_\alpha = 1$, and ignore the others. This eliminates any beneficial effect from a hit in one strip on hits in overlapping strips, so the two indicator variables are negatively correlated, $E(\bar{Y}_\alpha \bar{Y}_\beta) \leq \bar{p}_\alpha \bar{p}_\beta$, and

$$b_2 \leq L(2W-2)p_\alpha^2 = \lambda^2(2W-2)/L.$$

Combining this with the bound for $b_1$ in (6), we have

THEOREM 2. *Let $\bar{p}_\alpha = P(\bar{Y}_\alpha = 1)$ and $\bar{\lambda} = L\bar{p}_\alpha$. Let $\bar{V} = \sum_{\alpha \in I} \bar{Y}_\alpha$ be the number of occurrences under the new counting scheme and $Z$ be a Poisson random variable with $EZ = \bar{\lambda}$, then*

$$d_{\mathrm{TV}}(\bar{V}, Z) \leq \frac{2(4W - 3)}{L} \bar{\lambda}(1 - e^{-\bar{\lambda}}).$$

To compare with the results in the table above, when $L = 1024$ and $\bar{\lambda} = 1$, the bound is

$$0.02593 \text{ for } W = 6 \quad \text{and} \quad 0.03580 \text{ for } W = 8.$$

The good news about Theorem 2 is that the bound no longer depends on the word, and in case $W = 6$ is a dramatic improvement of the best case of the previous bound. The bad news is that it is difficult to analytically compute $\bar{\lambda}$. Based on results of Aldous [1], we guess that $\bar{\lambda} = \lambda/EC$, where $EC$ is the expected clump size. In Section 3.3 we derive an approximation for the clump size $E\hat{C}$ (where the hat indicates it is an approximation) by starting with the target word in the strip $A_\alpha$ for some $\alpha \in I$, random letters outside and computing the expected number of occurrences in overlapping strips $A_\beta$ for $\beta \in B_\alpha$ before the number of letters matching the target word hits 0 in that strip.

To evaluate the quality of the approximation, we turn to simulation. Since each letter changes at rate 1, the naive guess for the waiting time in one strip is $4^W/W$. If we had $L$ independent copies, this would reduce the time to $4^W/WL$. We will do our simulations for the embedded jump chain so this translates into $4^W/W$ mutations, which is 682.67 for $W = 6$ and 8192 for $W = 8$. Multiplying by our estimate of the clump size gives our prediction for the mean. In the table below, we compare with $E(T_A|T_A > 0)$ since the Poisson clumping heuristic estimates the mean time between occurrences. As simulation results given in the next table indicate, our formula is not very successful at predicting the mean number of mutations needed to produce the event $A = \{$the target word is found somewhere in the $L$ nucleotides$\}$ when $W = 6$, or when $W = 8$ and the word is $ACACACAC$. On the other hand, these results show that the naive estimate of $4^W/W$ mutations is never wrong by more than 20% (see Table 4).

To understand the reason for this, note that $E\hat{C}$ is computed under the assumption that when the word occurs the letters in adjacent positions are random. However, if we are looking at the first occurrence, then we are conditioning on the word having not been seen in adjacent strips.

TABLE 4

| Word | $E(T_A \mid T_A > 0)$ | $(4^W/W) \cdot E\hat{C}$ | $E\hat{C}$ |
|------|------|------|------|
| *AACCGT* | 717.32 | 770.97 | 1.129 |
| *ACGCTA* | 719.45 | 773.65 | 1.133 |
| *ACAGCA* | 729.49 | 785.50 | 1.150 |
| *AACGAA* | 732.79 | 797.97 | 1.171 |
| *AACAAC* | 746.42 | 823.96 | 1.210 |
| *ACACAC* | 806.85 | 900.34 | 1.318 |
| *ACAGCTGT* | 8674 | 8704 | 1.0624 |
| *ACAAGGGC* | 8685 | 8737 | 1.0665 |
| *ACAGACAG* | 8722 | 8881 | 1.0841 |
| *AAAACAAA* | 8825 | 8874 | 1.0832 |
| *AACAACAA* | 9013 | 9106 | 1.1116 |
| *ACACACAC* | 9584 | 10037 | 1.2253 |

c. *Waiting time distribution.* Putting together the results from parts a and b suggests:

APPROXIMATION 2b. *Under $P_\pi$ the waiting time for the target word to appear somewhere in one DNA sequence of length $L$ is*

$$\begin{cases} \approx (1 - e^{-\gamma})\delta_0 + e^{-\gamma}\xi, & \text{when } W = 6, \\ \approx \xi, & \text{when } W = 8, \end{cases}$$

*where $\gamma = L/4^W$ and $\xi$ an exponential distribution with mean $(4^W/LW)EC$.*

We have argued in the first part of this section that the size of the atom at 0 is accurately estimated by $1 - e^{-\gamma}$ when $W = 6$ and by 0 when $W = 8$. Consider first the case $W = 8$. Since the waiting time in one strip is approximately exponential with mean $4^W/W$ and the waiting time for the segment of length $L$ is roughly $4^W/LW$, if the time $T$ in the definition of $Y_\alpha$ is of this order of magnitude,

$$p_\alpha = P(Y_\alpha = 1) \approx TW/4^W,$$

where $W/4^W$ is the density of the exponential at 0. Taking into account the correction for clumping,

$$\bar{\lambda} = L\bar{p}_\alpha \approx Lp_\alpha/EC = TLW/(4^W EC).$$

When $W = 6$ this reasoning can be applied at positive times.

To check the prediction of an exponential distribution of positive values, Figures 1 and 2 show simulation results for waiting times for *AACCGT* and *ACACAC*. In each case, we group hitting times into bins of size 100 and plot the logarithm of the number of observations versus the hitting time count.

2.3. *W nucleotides in a population.*   Consider a population of $N$ diploid individuals. Following a standard practice in mathematical population genetics, we will formulate the dynamics of the Moran model as if there were $2N$ haploid individuals. Each individual, a string of $W$ letters from $\{A, C, G, T\}$, is replaced at rate 1, that is, individual $x$ lives for an exponentially distributed amount of time (with mean 1) and then is replaced. To replace individual $x$, we choose at random from the whole set of individuals, including $x$ itself, and make a copy. With probability $W\mu$, we randomly change one of the letters, and the new individual joins the population, replacing $x$.

We assume that, initially, everyone in the population has the same randomly chosen $W$ letter DNA word. To explain the reason for this assumption, we note that a standard result of population genetics implies that if we draw the genealogical tree tracing the entire population back to the most recent

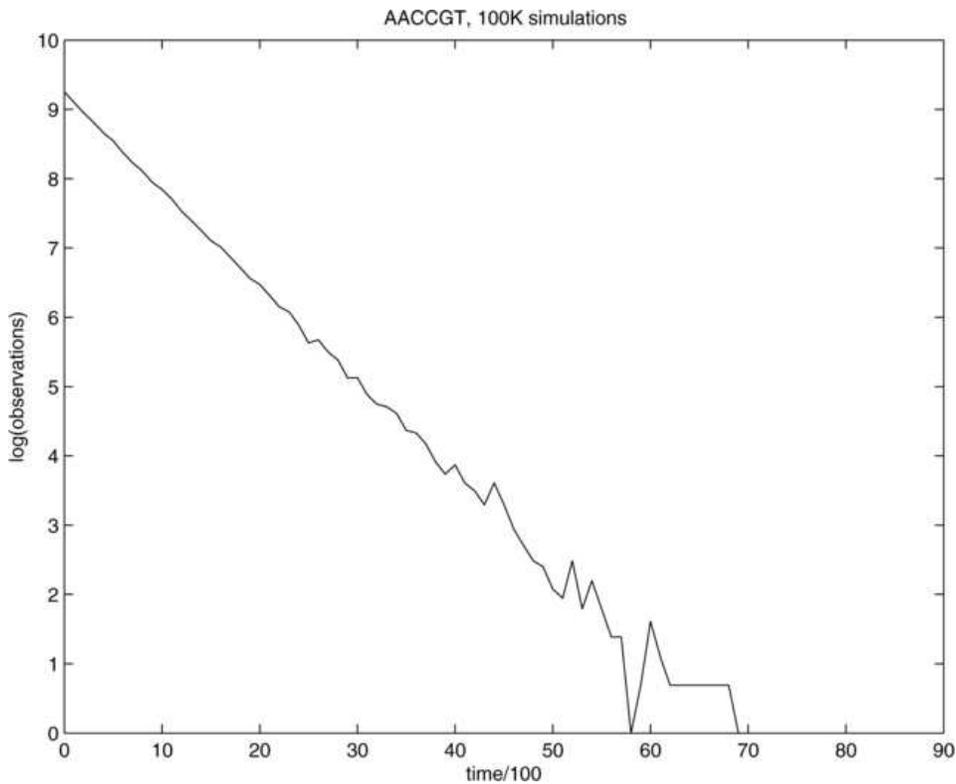

Fig. 1.   *Waiting time for target word AACCGT in* 100,000 *replications of the fixation chain simulation. Grouping hitting times into bins of size* 100, *we plot the logarithm of the number of observations versus the hitting time count. Despite the bad total variation bound, we see a good fit of the exponential distribution.*

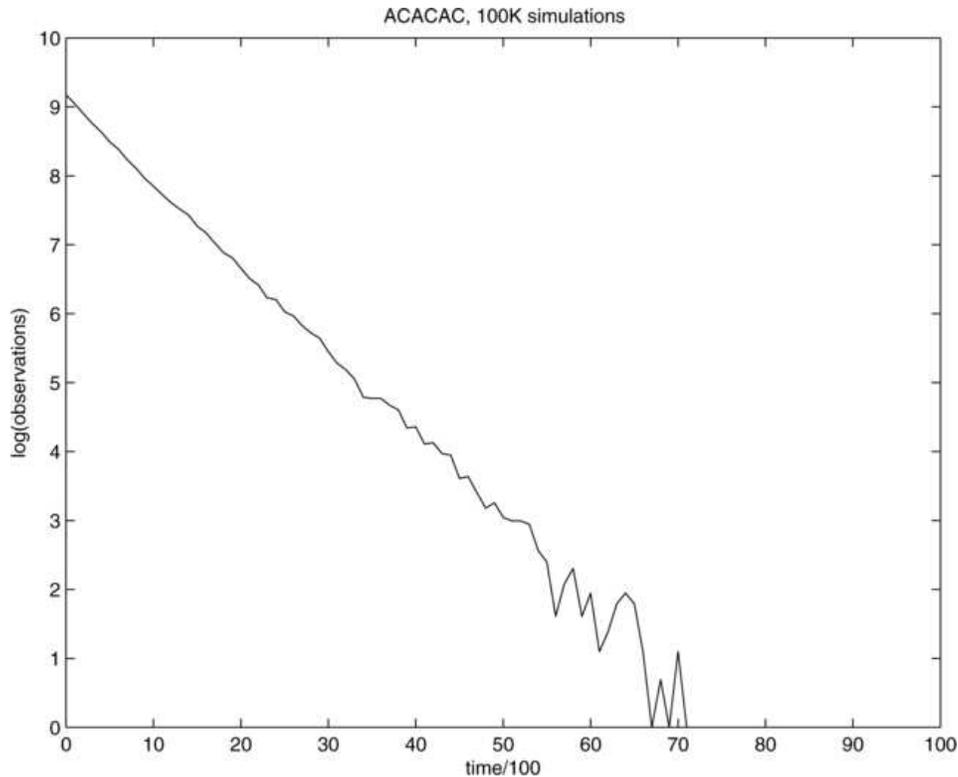

FIG. 2. *Waiting time for target word ACACAC in* 100,000 *replications of the fixation chain simulation. Grouping hitting times into bins of size* 100, *we plot the logarithm of the number of observations versus the hitting time count. Despite the bad total variation bound, we see a good fit of the exponential distribution.*

common ancestor, the expected total time in the tree is

$$(7) \qquad 2N \sum_{k=2}^{2N} k \cdot \frac{1}{\binom{k}{2}} = 4N \sum_{j=1}^{2N-1} 1/j \sim 4N \ln(2N).$$

The factor of $2N$ is the time scaling needed for a population of $N$ diploids to end up with the coalescent in which $k$ lineages coalesce to $k-1$ at rate $\binom{k}{2}$. For a more detailed explanation of this and the other facts from population genetics we make use of, the reader should consult Durrett [8], Ewens [11] or Tavaré [21].

From (7) we see that the expected number of mutations on the tree is $W\mu \cdot 4N \ln(2N)$. For reasons that will become clear in the next subsection, we will concentrate here on the case $W = 8$. Taking $\mu = 10^{-8}$ and $W = 8$, the result is 0.0316. In words, when $N = 10^4$, our value for the human population, 96.8% of the time there is no variation in the population.

Let $X_t(i) \in \{A, C, G, T\}^W$ be the state of individual $i$ at time $t$. Let $F = \{t \geq 0 : X_t(i) = X_t(1) \text{ for all } i\}$ be the set of times when all individuals have the same $W$ letter word. Here $F$ is for fixation, which is the genetics term for one word being present in all members of the population. Let $\mathcal{T}_n$ be the time of the $n$th fixation:

$$\mathcal{T}_n = \inf\{t > \mathcal{T}_{n-1} : t \in F, X_t(1) \neq X_{\mathcal{T}_{n-1}}(1)\},$$

where $\mathcal{T}_0 = 0$. Let $Y_n = X_{\mathcal{T}_n}(1)$. $Y_n$ is our discrete time fixation chain, which gives the state of the word in the population after the $n$th fixation. Let $L_n \in \{0, 1, \ldots, W\}$ be the number of letters in $Y_n$ that match the target word, and let $\tau_k = \inf\{n : L_n = k\}$ be the first time that $L_n$ hits state $k$. Let

$$\rho = \frac{2N\mu/9W}{1/2N + 2N\mu/9W} = \frac{4N^2\mu/9W}{1 + 4N^2\mu/9W}.$$

In terms of the fixation chain, we can state:

APPROXIMATION 3. *Let* $\xi_1, \xi_2, \ldots \in \{0, 1\}$ *be independent with* $P(\xi_i = 1) = \rho$. *Let* $S = \inf\{n : L_n = W - 1, \text{ or } L_n = W - 2 \text{ and } \xi_n = 1\}$. *The expected time to find the target word in* $W$ *nucleotides in a population of size* $N = 10^4$ *is* $\approx E_\pi S/(W\mu)$.

Using information on expected hitting times from the Appendix, we find

| $W$ | 6 | 8 |
|---|---|---|
| $E_\pi S$ | 214 | 2300. |

Since $\mu = 10^{-8}$, these numbers translate into huge waiting times: $3.567 \times 10^9$ and $2.875 \times 10^{10}$ generations, or 89.2 and 719 billion years respectively (using 25 years as the human generation time). This shows that it is important that regulatory sequences can occur in some region rather than at a fixed location.

REMARK. As calculations below will show, Approximation 3 becomes exact if $N \to \infty$ with $N^2\mu \to 1$ and is not accurate unless $N^3\mu^2$ is small. Thus, it is not valid for Drosophila populations with effective population size $N = 10^6$ and mutation rate $\mu = 10^{-8}$.

EXPLANATION. A new mutant is introduced into the population according to a Poisson process with rate $2N \cdot W\mu$, and it goes to fixation with probability $1/2N$. Thus, fixations occur at rate $W\mu$ and we will need of order $2N$ mutations before one reaches fixation. Since any $W$ letter word has only $3W$ one mutation neighbors, we can see that soon after $L_n = W - 1$ we will have a mutation that gives us the target word.

The next step is to consider the possibility that the target word will be reached when $L_n = W - 2$. Once a mutation occurs introducing a new word, the number of individuals with the mutant word evolves according to the Moran model. In Section 4 we will describe the dynamics of this model in detail and show that (for $N = 10^4$ and $\mu = 10^{-8}$), conditional on the mutation dying out, the probability of a second mutation before the first one dies out is $NW\mu = 8 \times 10^{-4}$. Taking into account that in order to increase $L_n$ from $W - 2$ to $W$, we need one of two possible mutations at the start and the right second mutation, the probability is

$$(8) \qquad \frac{2}{3W} \cdot NW\mu \cdot \frac{1}{3W} = \frac{2N\mu}{9W}.$$

Since fixation has probability $1/2N$, the probability of reaching the target word before the next fixation is $\rho$. To rule out success when $L_n = W - 3$, we use the reasoning for (8) to conclude that when $W = 8$, the expected number of good triple mutations (i.e., ones that will increase $L_n = W - 3$ to $W$) before we find the target word is $4.44 \times 10^{-4}$ (see Section 4 for more details).

At this point we have shown that if we measure time in units of the discrete time fixation chain, then Approximation 3 gives a very accurate approximation of the waiting time. To return to continuous time, we scale up by multiplying by the mean waiting time between fixations, and recall that, for any random variable, the mean of the sum is the sum of the means.

2.4. *A segment in a population.* We consider a population of $N$ diploid individuals, following the dynamics of the Moran model, but this time with a string of $L$ letters from $\{A, C, G, T\}$ indexed by the integers mod $L$ to avoid boundary effects. As in Section 2.2, we begin with:

a. *Almost matches in a random DNA sequence of length $L$.* Pick a sequence at random from $\{A, C, G, T\}^L$. In one window of width 8 the probability of matching 7 out of 8 letters is $8(3/4)(1/4)^7$, so the expected number of almost matches in $L = 1024 = 4^5$ nucleotides is $3/8$. Letting $M_i$ be the number of "matches minus $i$" in the initial condition, that is, words that disagree with the target word in $i$ letters, it is easy to see that

| $W$ | 6 | 8 |
|-----|------|--------|
| $EM_1$ | 4.5 | 0.375 |
| $EM_2$ | 33.75 | 3.9375. |

APPROXIMATION 4a.   *The number of matches minus $i$ in the initial configuration is approximately Poisson with mean $EM_i$, for $i = 1, 2$ as given in the table.*

In principle, one can use the methods in Section 3 to bound the error in these Poisson approximations, but the details get very messy and the bound depends in a complicated way on how the word overlaps with its shifts. For the 384 best words, the total variation distance is <0.01, and in more than 75% of cases the distance is <0.1, but in other cases one would need to use the Poisson clumping heuristic to get a good approximation.

b. *Quick success in a population.* The previous result was for one sequence, so to extend this to the population level, we use:

LEMMA 1. *When $W = 6$ the number of match minus 1's in the population at time 0 is $\approx 2NM$, where $M$ has a Poisson distribution with mean 4.5.*

Intuitively, the correlations between the various sequences are so strong that the number of copies in the population is with high probability roughly $2N$ times the number, $M$, in the most recent common ancestor of the population. We will give details in Section 5. A corollary of this computation is that if we define the population consensus sequence to be the most common nucleotide at each position, then with high probability the population consensus coincides with the sequence of the most recent common ancestor.

Consider first the situation when $M = 1$, so that there are roughly $2N$ match minus 1's in the population. When $N = 10^4$ and $\mu = 10^{-8}$, the rate for getting the right mutation at the desired location is $2N\mu/3 = (2/3) \times 10^{-4}$, so the expected time is roughly exponential with mean 15,000 generations or 375,000 years. If $M = k$, then the waiting time is the minimum of $k$ independent copies of the waiting time when $M = 1$, so the waiting time is divided by $k$. When $W = 6$ the number of copies is roughly Poisson with mean 4.5, so the probability of at least 1 is 0.988. Ignoring the zero cell, the expected waiting time is approximately

$$(9) \qquad 375{,}000 \cdot \sum_{k=1}^{\infty} e^{-4.5} \frac{4.5^k}{k!} \cdot \frac{1}{k} = 107{,}697 \text{ years.}$$

Thus, words of length 6 have an average waiting time of about 100,000 years, as quoted in the abstract. Note that the waiting time is not exponential, but is a mixture of exponentials.

c. *When there is no almost match.* Having taken care of the case $W = 6$, we now focus on $W = 8$. Approximation 3 says that waiting for the word to appear at a specified $W$ nucleotides in a population of size $N$ is almost the same as waiting for the death of a randomly killed version of the fixation chain. Thus, to estimate the waiting time where there is not a match minus 1 in the initial condition, we consider the fixation chain [defined as in Section 2.3, but this time for $X_t(i) \in \{A, C, G, T\}^L$] that jumps when the population is fixed for a new nucleotide at some position. This chain jumps at rate $L\mu$, so the probabilities of reaching the target word before the next fixation change from the values we have computed previously. When there is only one mismatch in a window, the target word appears in some individual at rate $2N\mu/3$, so the probability this occurs before the next fixation is

$$\rho_1 = \frac{2N\mu/3}{2N\mu/3 + L\mu} = (1 + 3L/2N)^{-1} = 20/23$$

$$\text{when } L = 1000 \text{ and } N = 10{,}000.$$

Now consider having two mismatches in a window. Since mutations occur at rate $2N\mu$ per nucleotide, those that fix one of the two wrong letters occur at rate $2N\mu \cdot 2/3$. The probability of correcting the other letter before the mutation dies out is $N\mu \cdot 1/3$, so the probability the target word appears in some individual before the next fixation is

$$\rho_2 = \frac{4(N\mu)^2/9}{4(N\mu)^2/9 + L\mu} = \frac{r}{r+1},$$

where

$$r = \frac{4N^2\mu}{9L}.$$

When $L = 1000$ and $N = 10{,}000$, $\rho_2 = 4/9000$.

$\rho_1 \approx 1$ and in most cases the next few fixations will not affect the word with seven matches, so we will stop with probability 1 when a match minus 1 is found. Let $T_D$ be the death time of the killed chain in which each match minus 2 (even those that overlap) independently results in the end of the process with probability $\rho_2 = 4/9000$. The choice in parentheses makes programming the simulation easier, but leads to a slight underestimate of the actual ending time.

APPROXIMATION 4b.   *The expected time to find the target word in a population of size $N = 10^4$ when each individual has $L$ letters is $\approx E_\pi T_D/(L\mu)$.*

To evaluate $E_\pi T_D$, we turn to simulation. The following results in Table 5 are based on 100,000 replications of the killed fixation chain.

To interpret Table 5, recall that the predicted probability of at least one match minus 1 in the initial configuration of the $L$ letter fixation chain is $1 - e^{-3/8} = 0.3127$, which is $\approx 5/16 = 0.3125$. Thus, there is no match minus 1 roughly 11/16's of the time. The smallest mean in the table is about 260. When $L = 1000$ and $\mu = 10^{-8}$, new fixations happen at rate

Table 5

| Word | $P_\pi(T_D = 0)$ | $E_\pi(T_D\|T_D > 0)$ |
|------|------------------|----------------------|
| $ACAGCTGT$ | 0.3199 | 259.95 |
| $ACAAGGGC$ | 0.3225 | 260.51 |
| $ACAGACAG$ | 0.3174 | 275.54 |
| $AAAACAAA$ | 0.3116 | 297.62 |
| $AACAACAA$ | 0.3030 | 293.93 |
| $ACACACAC$ | 0.2744 | 329.70 |

$L\mu = 10^{-5}$, so this translates into $2.6 \times 10^7$ generations or 650 million years. When there is exactly one match minus 1 in the initial condition, the waiting time is exponential with mean 375,000 years. Hence, the waiting time for words of length 8 follows the mixed distribution quoted in the abstract.

Figures 3 and 4 show simulation results for waiting times for $ACAGCTGT$ and $ACACACAC$ when there is no match minus 1 in the initial configuration, with hitting times grouped into bins of size 10. The simulation shows an exponential distribution for $(T_D|T_D > 0)$ under $P_\pi$, even though we do not have a good reason why the waiting times should have the lack of memory property. For the readers who are disappointed to see this final result done by simulation, we would like to note that if instead of considering the killed fixation chain, we did a simulation of the Moran model for 10,000 diploids for $2.6 \times 10^7$ generations, this would take $5.2 \times 10^{11}$ simulation steps. Thus, even one simulation would be time consuming and 100,000 unthinkable.

2.5. *Conclusions.* The calculations in Section 2.4 give the answer to Stone and Wray's [20] question: How long do we have to wait for an exact match to a given $W$ letter word to appear in a segment of length $L$ in some individual in a population of $N$ diploids? When $W = 6$ the mean waiting time is about 100,000 years, which easily allows differences to accumulate in the 6 million years between the divergence of humans and chimpanzees. For $W = 8$, we note that, in combination with Lemma 1, our results have shown that the mean waiting time is about 375,000 years if there is a match minus 1 in the population consensus sequence, but 650,000,000 years if there is not.

The second answer says that unless you are lucky enough to have a match minus 1 in the population consensus sequence, regulatory sequence changes can take an extremely long time. However, waiting for exact matches is not the right question for the evolution of regulatory sequences. The DNA sequence does not have to be exactly right for transcription factor binding to occur. The binding energy between a transcription factor and its DNA binding site is, to a good approximation, the sum of independent contributions from a small number of important positions in the binding site sequence $\sum_i \varepsilon_i$; see [12]. In a commonly used approximation called the two state model $\varepsilon_i \in \{0, \varepsilon\}$ (see [5]), the binding energy is determined by the number of mismatches $r$ between the two strings. The binding probability is commonly modeled by the Fermi function

$$p = \frac{1}{1 + \exp(\varepsilon(r - r_0))},$$

where $r_0$ is the threshold for a binding probability of $1/2$, and contrary to the usual mathematical usage $\varepsilon \approx 2$.

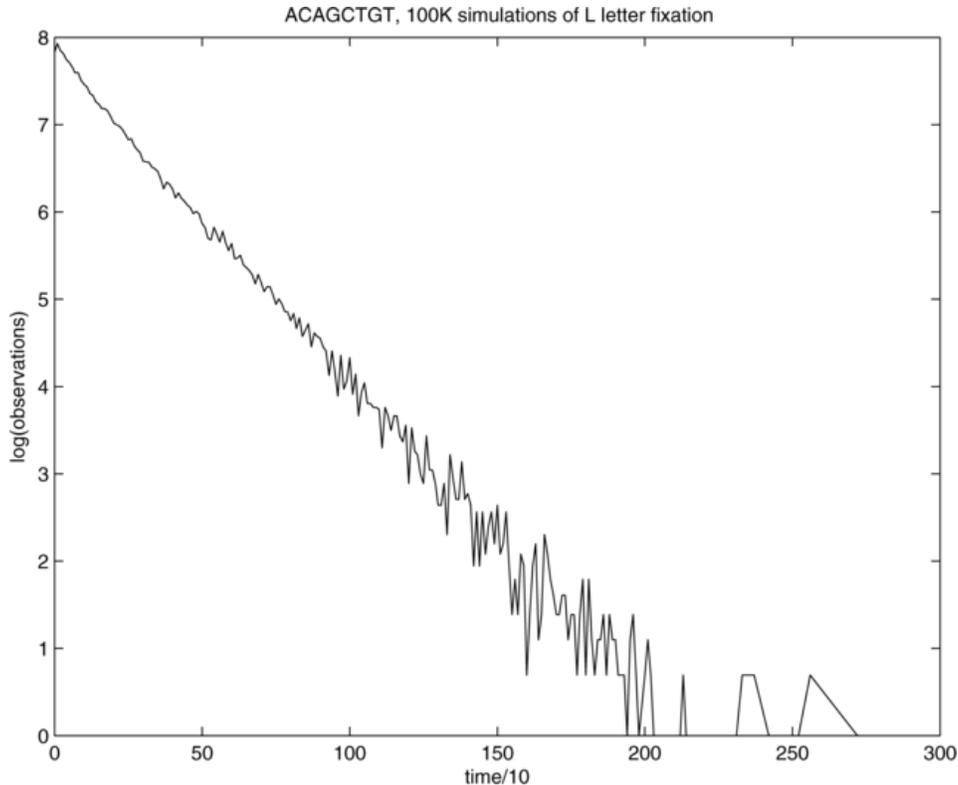

FIG. 3. *Waiting time for target word ACAGCTGT in* 100,000 *replications of the L letter fixation chain. We group hitting times into bins of size* 10 *and plot the logarithm of the number of observations versus the hitting time count. Again, we see a good fit of the exponential distribution.*

Based on this analysis, it seems reasonable to simplify the biological complexities of binding of transcription factors by saying that a match minus 1 is adequate. If we do this, then it follows from the analysis above that it

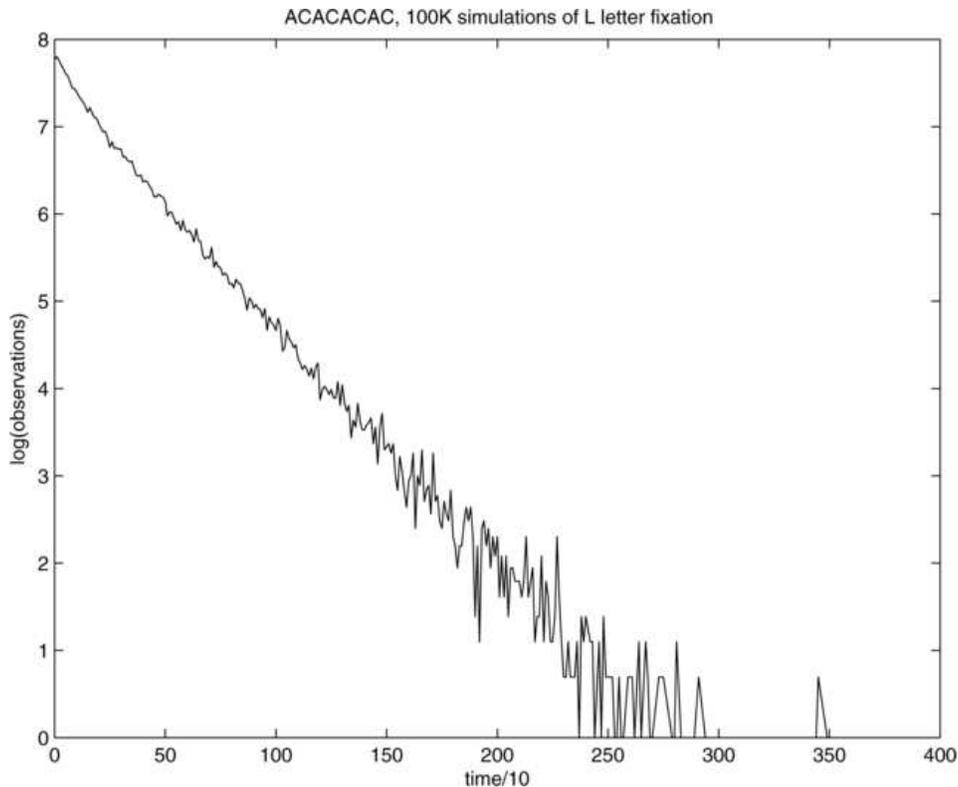

Fɪɢ. 4.   *Waiting time for target word ACACACAC in* 100,000 *replications of the L letter fixation chain. We group hitting times into bins of size* 10 *and plot the logarithm of the number of observations versus the hitting time count. Again, we see a good fit of the exponential distribution.*

is enough to have a match minus 2 in the fixation chain. From Approximation [4a], the mean number of these in the initial distribution is 33.75 when $W = 6$, and 3.94 when $W = 8$, so the very long waiting times are unlikely. In the case of an 8 letter word, generalizing [9] shows that the waiting time has mean

$$\frac{375,000}{2} \cdot \sum_{k=1}^{\infty} e^{-3.94} \frac{3.94^k}{k!} \cdot \frac{1}{k}$$

(10)

$$= 61,560 \text{ years.}$$

Here, we have divided by 2 since there are two sites where a mutation can upgrade a 6 out of 8 match to 7 out of 8.

Having announced our results, we should mention two things, each of which could change the answer by a factor of 100, but in different directions.

First, we have assumed that once the target word occurs in an individual, it sweeps to fixation in the population, that is, not long after introduction the frequency will rise to 1. In reality the probability of fixation is approximately the selective advantage conferred by the mutation $s$ and even for strongly beneficial mutations we have $s \leq 0.01$. This means that the mutation would need to arise more than 100 times in order to achieve fixation, which would increase the waiting time to 6 million years. In the other direction, our study has focused on changes in the regulation of one particular gene, but there are more than 20,000 genes in humans and chimpanzees, and changes in even 1% of these genes could be enough to explain the observed differences.

We have pursued an approach in which mutations are neutral, that is, they do not change the fitness of the individual. Several researchers have formulated models for the evolution of regulatory sequences, making use of models of the binding of regulatory proteins to DNA, similar to the ones described above, and assuming the fitness is proportional to the binding probability. Gerland and Hwa [13] assumed an infinite population and developed a theory based on the quasispecies approach of [10]. MacArthur and Brookfield [17] and Berg, William and Lassig [4] have developed related models that incorporate genetic drift. These models use different mechanisms, but also conclude that the existence of presites (which can be converted into binding sites by one mutation, i.e., our match minus 1's) largely determine where binding sites will evolve.

The results in this paper depend heavily on the fact that the population size $N = 10^4$. When $N \approx 10^5$ the estimates in Section 2.3 show that we cannot ignore the possibility of triple mutations between fixations. For $N = 10^6$ which is appropriate for Drosophila, the collection of sequences in the population is even more fluid. Experimental results for the even-skipped stripe 2 enhancer (see [16] and references therein) suggest that other evolutionary mechanisms are at work in this case. Carter and Wagner [6] have shown that in this case binding domains can shift through a combination of a mildly deleterious mutation followed by subsequent selection for a compensatory mutation. Prudhomme et al. [18] have shown recently that regulatory changes have caused two independent gains and five losses of wing pigmentation spots in the *Drosophila melanogaster* species group. Understanding the mechanisms at work in the evolution of gene regulation when $N = 10^6$ or for larger population sizes appropriate for bacterial or viral genomes is an interesting open problem.

**3. Chen–Stein calculations.** To fill in the missing details in our treatment of a segment of length $L$ in one DNA sequence (Section 2.2), the first detail is to consider the following:

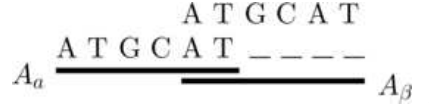

3.1. *Success in the initial condition.* For $\alpha, \beta \in I$ and $A_\alpha$, $A_\beta$ defined as in Section 2.2, when $\beta = \alpha + k$ for $1 \leq k \leq W - 1$, $A_\alpha$ and $A_\beta$ overlap in $W - k$ letters. Figure 5 shows an example with $W = 6$ and $k = 4$.

Let $y_k = 1$ if the word and the shifted word match exactly on the overlap, 0 otherwise. It is easy to see that $E[Y_\alpha Y_\beta] = y_k 4^{-(W+k)}$, so

$$b_2 = \sum_\alpha \sum_{\beta \in B_\alpha \setminus \{\alpha\}} E[Y_\alpha Y_\beta]$$

$$= \sum_\alpha 4^{-W} \cdot 2 \sum_{k=1}^{W-1} y_k 4^{-k}$$

$$= 2\gamma Q_0 \qquad \text{where } \gamma = L4^{-W} \text{ and } Q_0 = \sum_{k=1}^{W-1} y_k 4^{-k}.$$

$Q_0$ is the expected number of overlapping matches when shifting to the right conditioned on finding a match in the first $W$ nucleotides. The 2 comes from the symmetric case of shifting left as well as right. From this formula, it is straightforward to calculate the values in Table 2 in Section 2.2.a.

3.2. *Number of occurrences by time $T$.* Let $p_{\alpha\beta} = E[Y_\alpha Y_\beta]$. We estimate

$$p_{\alpha\beta} = p_{\alpha=\beta} + p_{\alpha<\beta} + p_{\beta<\alpha},$$

where $p_{\alpha=\beta}$ refers to hitting the target word at the same time in two overlapping strips $A_\alpha$ and $A_\beta$, while $p_{\alpha<\beta}$ (which is equivalent to the case $p_{\beta<\alpha}$) refers to hitting the target word in $A_\alpha$ strictly before $A_\beta$.

Consider the target word and its shift by $k$ letters, as in the picture above the definition of $y_k$. Define

$m_k = \#$ of matching letters inside the overlap of $W - k$ letters.

Consider the moment that we first get a match in $A_\alpha$ and let

$R_k = \#$ of letters matching the target word in $A_\beta \setminus A_\alpha$.

Since we start from a random initial sequence, $R_k$ is Binomial$(k, 1/4)$.

We start with $p_{\alpha=\beta}$. Let $y_k = 1$ if $m_k = W - k$, and 0 otherwise

$$\sum_{\beta \in B_\alpha \setminus \{\alpha\}} p_{\alpha=\beta} = 2p_\alpha Q, \qquad \text{where } Q = \sum_{k=1}^{W-1} y_k \left(\frac{1}{4}\right)^k \cdot \frac{W-k}{W},$$

where $(W - k)/W$ is the probability the last mutation was in the overlap, a necessary condition for hits in $A_\alpha$ and $A_\beta$ to occur at the same time. From this, it follows that

$$(11) \qquad \sum_\alpha \sum_{\beta \in B_\alpha \backslash \{\alpha\}} p_{\alpha = \beta} = 2\lambda Q.$$

Turning to the case $p_{\alpha < \beta} = p_{\beta < \alpha}$, let $r(k, j) = P(R_k = j)$, let $h(x) = P_x(T_W < T_0)$ and

$$z(m_k, k) = \sum_{j=0}^k r(k, j) h(m_k + j) \mathbb{1}_{\{m_k + j < W\}}.$$

This is the probability of hitting the target word in $A_\beta$ before the number of matches in $A_\beta$ returns to 0, when there are $m_k$ matches inside the overlap. We eliminate matches at time 0 since they are counted in $p_{\alpha = \beta}$. Since after the return to 0 there are $\leq T$ units of time left, $p_{\alpha + k}$ is an upper bound on the probability of a match after the return to 0 matches and before time $T$,

$$\sum_{\beta \in B_\alpha \backslash \{\alpha\}} p_{\alpha < \beta} \leq 2p_\alpha \sum_{k=1}^{W-1} [z(m_k, k) + p_{\alpha + k}].$$

Letting $S = \sum_{k=1}^{W-1} z(m_k, k)$ and noting that $2 \sum_\alpha p_\alpha \sum_{k=1}^{W-1} p_{\alpha+k}$ is essentially $b_1$ except we get $2W - 2$ instead of $2W - 1$, we have

$$(12) \qquad \sum_\alpha \sum_{\beta \in B_\alpha \backslash \{\alpha\}} [p_{\alpha < \beta} + p_{\beta < \alpha}] \leq 2 \left( \frac{(2W - 2)\lambda^2}{L} + 2\lambda S \right).$$

Combining (11) and (12), we have

$$b_2 \leq \lambda^2 (4W - 4)/L + \lambda(2Q + 4S).$$

Using this, we can compute the values in Table 3.

3.3. *Clump size.* Our next topic is to explain our approximation for the expected clump size:

$$(13) \qquad E\hat{C} \approx \frac{1 + 2(S + Q)}{1 - a}.$$

Intuitively, the numerator gives the expected number of strips in which we get a hit before the number of matching letters in that strip returns to 0. From the definition in the previous subsection, we see that $2Q$ takes care of the hits that occur at the same time, $2S$ the ones that come later. Multiplying by $1/(1 - a)$ accounts for occurrences after the first one in a strip.

**4. Population dynamics.** Here we are concerned with the derivation of Approximation 3. In that setting we have a population of $N$ diploid individuals, each of which has a $W$ letter DNA word. Suppose that at time 0 all individuals in the population have the same $W$ letter word except for one mutant that differs in one letter, and for the moment ignore the possibility of additional mutations. Let $Z_t$ = the number of "mutants" in the population at time $t$, that is, the number of individuals with the mutated word. $Z_t$ follows the dynamics of the Moran model. That is, when $Z_t = k$,

    I. replace a mutant with a mutant at rate $k \cdot \frac{k}{2N}$,
    II. replace a mutant with a nonmutant at rate $k \cdot \frac{2N-k}{2N}$,
    III. replace a nonmutant with a mutant at rate $(2N - k) \cdot \frac{k}{2N}$,
    IV. replace a nonmutant with a nonmutant at rate $(2N - k) \cdot \frac{2N-k}{2N}$.

Let $T_k = \inf\{t \geq 0 : Z_t = k\}$ be the hitting time of $k$. We stop the process when $Z_t$ hits 0 (the mutation dies out) or $2N$ (the mutation becomes fixed in the population) and call the result an "excursion." Our first goal in this section is to compute the expected value of $B$ = the total number of mutant births in one excursion, that is, the number of times I or III occurs before the process enters one of the absorbing states 0 or $2N$. Writing $a_N \sim b_N$ to indicate $a_N/b_N \to 1$ as $N \to \infty$, we have

LEMMA 2. $E_1(B|T_0 < T_{2N}) \sim N$.

PROOF. To prove this claim, we begin by noting that up and down jumps occur at the same rate, so the embedded jump chain (which gives the sequence of states visited by the continuous time chain) is a simple random walk $S_n$ with $S_0 = 1$.

*Simple random walk.* Again let $T_k = \inf\{n \geq 0 : S_n = k\}$ be the hitting time of $k$. Let $N_k$ = number of visits to state $k$ for $S_n$ before $T_0$. Our first step is to compute

$$E_1(N_k|T_0 < T_{2N}) = \frac{P_1(T_k < T_0|T_0 < T_{2N})}{P_k(T_k^+ > T_0|T_0 < T_{2N})},$$

where $T_k^+ = \inf\{n \geq 1 : S_n = k\}$ is the return time. To explain the formula, the numerator gives the probability the mutants achieve a population of size $k$. Once this occurs, the simple random walk has a geometric number of visits to $k$ with success probability $P_k(T_k^+ > T_0|T_0 < T_{2N})$. Using the definition of conditional probability and the Markov property,

$$P_1(T_k < T_0|T_0 < T_{2N}) = \frac{P_1(T_k < T_0, T_0 < T_{2N})}{P_1(T_0 < T_{2N})}$$

$$= \frac{P_1(T_k < T_0)P_k(T_0 < T_{2N})}{P_1(T_0 < T_{2N})}$$

(14)
$$= \frac{(1/k)(1 - k/(2N))}{1 - 1/(2N)}$$

$$= \frac{2N - k}{k(2N - 1)}.$$

The next to last equality comes from the fact that $S_n$ is a martingale and, hence, $P_i(T_j < T_0) = i/j$. Turning to the denominator of $E_1(N_k|T_0 < T_{2N})$,

(15)
$$P_k(T_k^+ > T_0|T_0 < T_{2N}) = \frac{P_k(T_k^+ > T_0, T_0 < T_{2N})}{P_k(T_0 < T_{2N})}$$

$$= \frac{(1/2) \cdot (1/k)}{1 - k/(2N)} = \frac{2N}{2k(2N - k)}.$$

The $1/2$ comes from the fact that $S_n$ has to jump from $k$ to $k - 1$ in order to be able to hit $0$ before going back to $k$. Combining our results gives

(16)    $E_1(N_k|T_0 < T_{2N}) = \left( \dfrac{2N - k}{k(2N - 1)} \right) \left( \dfrac{2k(2N - k)}{2N} \right) = \dfrac{2(2N - k)^2}{2N(2N - 1)}.$

*Moran model.* Next we want to compute the mean number of mutant births while the Moran model is in state $k$. Let $W_k$ be the number of type I events before a jump makes the chain leave $k$. Let $J_k^i$ be the number of type $i$ jumps while the chain is in state $k$:

$$E_1(J_k^1|T_0 < T_{2N}) = E(N_k|T_0 < T_{2N}) \cdot EW_k.$$

$W_k$ has a shifted geometric distribution with mean $EW_k = 1/P(\text{jump}) - 1$, where

$$P(\text{jump}) = \frac{2k(2N - k)/(2N)}{kk/(2N) + 2k(2N - k)/(2N)} = \frac{4N - 2k}{4N - k},$$

which implies $EW_k = k/(4N - 2k)$. Using (16) and summing, we have

$$\sum_{k=1}^{2N-1} E_1(J_k^1|T_0 < T_{2N}) = \sum_{k=1}^{2N-1} \frac{2(2N - k)^2}{2N(2N - 1)} \cdot \frac{k}{4N - 2k}$$

$$= 2N \sum_{k=1}^{2N-1} \left( 1 - \frac{k}{2N} \right) \cdot \frac{k}{2N - 1} \cdot \frac{1}{2N}$$

$$\sim 2N \int_0^1 (1 - x)x \, dx = N/3.$$

To compute the expected number of type III jumps in state $k$, we begin by noting that the probability we ever visit $k$ is given by (14). Using the calculation for (15),

$$P_k(\text{no up jump from } k \text{ before } T_0 | T_0 < T_{2N})$$

(17)
$$= P_k(T_{k+1} > T_0 | T_0 < T_{2N}) = \frac{1/(k+1)}{1 - k/(2N)} = \frac{2N}{(k+1)(2N-k)}.$$

Since the number of up jumps from $k$ for the conditioned chain has a shifted geometric distribution, using (14) gives

$$E_1(J_k^3 | T_0 < T_{2N}) = \frac{2N-k}{k(2N-1)} \left( \frac{(k+1)(2N-k)}{2N} - 1 \right).$$

If $k$ and $2N - k$ are large, then $k(2N-k)/2N$ is much bigger than 1 and the $-1$ can be ignored. Summing, we have

$$\sum_{k=1}^{2N-1} E_1(J_k^3 | T_0 < T_{2N}) \sim 2N \sum_{k=1}^{2N-1} \frac{(2N-k)^2}{2N(2N-1)} \cdot \frac{1}{2N}$$

$$\sim 2N \int_0^1 (1-x)^2 \, dx = 2N/3.$$

Adding this to the expected number of jumps of type I gives the desired result. $\square$

We next consider the final excursion in which the new mutant fixes.

LEMMA 3. $E_1(B | T_{2N} < T_0) \sim 2N^2$.

PROOF. We use the notation of the previous proof. Here, $P_1(T_k < T_0 | T_{2N} < T_0) = 1$, so

$$E_1(N_k | T_{2N} < T_0) = \frac{1}{P_k(T_k^+ > T_{2N} | T_{2N} < T_0)}.$$

The denominator

$$= \frac{(1/2)P_{k+1}(T_{2N} < T_k)}{P_k(T_{2N} < T_0)} = \frac{(1/2) \cdot 1/(2N-k)}{k/(2N)} = \frac{2N}{2k(2N-k)},$$

since for simple random walk $P_{k+1}(T_{2N} < T_k) = P_1(T_{2N-k} < T_0)$. To compute the expected number of jumps of type I in state $k$, we note that

$$E_1(J_k^1 | T_{2N} < T_0) = E(N_k | T_{2N} < T_0) \cdot EW_k = \frac{2k(2N-k)}{2N} \cdot \frac{k}{4N-2k} = \frac{k^2}{2N}.$$

Summing, we have

$$\sum_{k=1}^{2N-1} E_1(J_k^1 | T_{2N} < T_0) \sim (2N)^2/3.$$

To compute the expected number of jumps of type III in state $k$, we begin by using the computation in (17) to conclude

$$P_{k+1}(\text{no down jump from } k+1 \text{ before } T_{2N}|T_{2N} < T_0)$$
$$= P_{k+1}(T_k > T_{2N}|T_{2N} < T_0)$$
$$= P_{2N-k-1}(T_{2N-k} > T_0|T_0 < T_{2N}) = \frac{2N}{(k+1)(2N-k)}.$$

Thus, the expected number of down jumps from $k+1$ (which come from type II changes):

$$E_1(J_{k+1}^2|T_{2N} < T_0) = \frac{(k+1)(2N-k)}{2N} - 1.$$

Since on $\{T_{2N} < T_0\}$ this is one less than the expected number of up jumps from $k$,

$$E_1(J_k^3|T_{2N} < T_0) = \frac{(k+1)(2N-k)}{2N}.$$

Summing, we have

$$\sum_{k=1}^{2N-1} E_1(J_k^3|T_{2N} < T_0) \sim (2N)^2 \int_0^1 x(1-x)\,dx = (2N)^2/6.$$

Adding this to the expected number of jumps of type I gives the desired result.   $\square$

*Derivation of Approximation* 3. As explained in Section 2.3, fixations occur at rate $W\mu$ and soon after $L_n = W - 1$ we will have a mutation that gives us the target word.

The next step is to consider the possibility that the target word will be reached when the fixation chain $L_n = W - 2$. Lemma 2 shows that if $B$ is the number of times a new mutant is born in an excursion, then

$$E_1(B|T_0 < T_{2N}) \sim N.$$

Given this result, if $L_n = W - 2$, the expected number of times that we will have a correct mutation at the start and the right second mutation before the first one dies out is

$$\frac{2}{3W} \cdot N \cdot W\mu \cdot \frac{1}{3W} = \frac{2N\mu}{9W}.$$

The first factor gives the fraction of the $3W$ possible mutations that fix one of the mismatches, the second, the expected number of births per excursion that can result in mutation (by Lemma 2), the third, the mutation probability and the fourth factor the fraction of the $3W$ possible mutations that fix the one remaining mismatch.

Since fixation has probability $1/2N$, thinking about a race between two rare events shows that the probability of reaching the target word before the next fixation is approximately

$$\rho = \frac{2N\mu/9W}{1/2N + 2N\mu/9W} = \frac{4N^2\mu/9W}{1 + 4N^2\mu/9W}.$$

To verify this, we note that the probability of not reaching the target word is

$$
\begin{aligned}
(18) \qquad & \sum_{k=0}^{\infty} \left(1 - \frac{2N\mu}{9W}\right)^k (1/2N)(1 - 1/2N)^k \\
& \approx \sum_{k=0}^{\infty} (1/2N)\left(1 - \frac{2N\mu}{9W} - \frac{1}{2N}\right)^k \\
& = \frac{1/(2N)}{2N\mu/(9W) + 1/(2N)} = \frac{1}{1 + 4N^2\mu/9W}.
\end{aligned}
$$

An important consequence of this calculation is that we can ignore the interaction between mutations. When $W = 8$ and $\mu = 10^{-8}$, the probability that between two successive fixations a second mutation arises before the previous one dies out is $1/19$. The probability this mutation reaches size $K$ in the population is $1/K$, so with high probability none will reach a large size. There is also the issue of second mutations while a successful mutation is reaching fixation. As Lemma 3 shows,

$$E_1(B|T_{2N} < T_0) \sim 2N^2.$$

Thus, the expected number of second mutations during this excursion is $2N^2\mu = 2$ (when $N = 10^4$ and $\mu = 10^{-8}$), but with high probability none will reach a large size. We ignore the possibility of finding the target word in the final excursion when $L_n = W - 2$, since this can happen only if $L_{n+1} = W - 1$, in which case the killed fixation chain will terminate at time $n + 1$.

To rule out success when $L_n = W - 3$, we extend the reasoning used for $L_n = W - 2$ to conclude that the expected number of good triple mutations (i.e., ones that will increase $L_n = W - 3$ to $W$, bringing us to the target word) between two fixations is

$$2N \cdot \frac{3}{3W} \cdot N \cdot W\mu \cdot \frac{2}{3W} \cdot N \cdot W\mu \cdot \frac{1}{3W} = \frac{4N^3\mu^2}{9W}.$$

The first factor is the expected number of excursions between two fixations, the second the fraction of the $3W$ possible mutations that fix one of the mismatches, the third the expected number of births per excursion that can result in mutation, the fourth the mutation probability, the fifth the fraction of the $3W$ possible mutations that fix one of the two remaining mismatches and so on. When $W = 8$, using the fact established in (25) that $E(N_5(\tau_7)) = 80$, we see that the probability of a good triple mutation before we find the target word is $4.44 \times 10^{-4}$.

**5. Proof of Lemma 1.** The calculation for (7) has shown that the expected number of mutations in a segment of length $L$ in a population of $N$ diploids is $L\mu \cdot 4N \ln(2N) = 3.95$ when $L = 1000$, $N = 10{,}000$ and $\mu = 10^{-8}$. Since 4.5 windows of length 6 contain 27 nucleotides, this implies that the expected number of match minus 1's that are disrupted by mutation is

$$27 \cdot \frac{3.95}{1000} = 0.1066.$$

In the other direction, since there are an average of 33.75 match minus 2's and each has only two sites that can be fixed by one of three possible mutations, the number of match minus 1's created is

$$33.75 \cdot 2 \cdot \frac{1}{3} \cdot \frac{3.95}{1000} = 0.0888.$$

The last two expected values are not negligible, but mutations when they exist are rare. Well-known results about the site frequency spectrum imply that the expected number of individuals with the mutant nucleotide at a site given that a mutation has occurred is

$$\sum_{k=1}^{2N} \frac{1}{k \ln(2N)} k = \frac{2N}{\ln(2N)} = 0.1009(2N)$$

when $N = 10{,}000$.

## APPENDIX: COMPUTATIONS FOR THE MUTATION CHAIN

Our first task is to explain the computation of the mixing time $\tau_2$ that appears in the derivation of Approximation 1. Let

$$U = \begin{pmatrix} 1/4 & 1/4 & 1/4 & 1/4 \\ 1/4 & 1/4 & 1/4 & 1/4 \\ 1/4 & 1/4 & 1/4 & 1/4 \\ 1/4 & 1/4 & 1/4 & 1/4 \end{pmatrix}, \qquad V = \begin{pmatrix} 0 & 1/3 & 1/3 & 1/3 \\ 1/3 & 0 & 1/3 & 1/3 \\ 1/3 & 1/3 & 0 & 1/3 \\ 1/3 & 1/3 & 1/3 & 0 \end{pmatrix}.$$

The eigenvalues of $U$ are 1, 0, 0, 0, since any vector $x$ orthogonal to the constant vector has $Ux = 0$. To compute the eigenvalues of $V$, we note that $V = (4/3)U - (1/3)I$, so its eigenvalues are 1, $-1/3$, $-1/3$, $-1/3$ and the spectral gap is $4/3$. If we consider the discrete time chain in which there are $W$ letters and one changes at each step, the spectral gap is $4/3W$ or $\tau_2 = 3W/4$. Of course, if we speed up the chain so that each coordinate jumps at rate 1, then $\tau_2 = 3/4$.

Let $X_n$ have state space $S = \{0, 1, \ldots, W\}$ and transition probabilities

(19)
$$
\begin{aligned}
p(x, x-1) &= x/W, \\
p(x, x+1) &= (1/3)(W-x)/W, \\
p(x, x) &= (2/3)(W-x)/W,
\end{aligned}
$$

where all other $p(x, y) = 0$. In this section we will use standard Markov chain techniques to compute several quantities of interest for this Markov chain.

Let $h(x) = P_x(T_W < T_0)$ be the probability of hitting $W$ before hitting state 0 when starting in state $x$. $h(x)$ satisfies

(20) $\quad h(x) = p(x, x+1)h(x+1) + p(x, x-1)h(x-1) + p(x, x)h(x)$

for $0 < x < W$, with boundary conditions $h(W) = 1$ and $h(0) = 0$. We solve for $h(x)$ recursively. Since $1 - p(x, x) = p(x, x+1) + p(x, x-1)$, rearranging gives

$$
p(x, x-1)[h(x) - h(x-1)] = p(x, x+1)[h(x+1) - h(x)],
$$

which implies

(21)
$$
\begin{aligned}
h(x) - h(x-1) &= \left[\frac{p(x, x+1)}{p(x, x-1)}\right][h(x+1) - h(x)] \\
&= \left[\frac{W-x}{3x}\right][h(x+1) - h(x)].
\end{aligned}
$$

Setting $h(W) - h(W-1) = C$ and using the recursion in (21) gives

$$
h(W-1) - h(W-2) = \left[\frac{1}{3(W-1)}\right]C,
$$

$$
h(W-2) - h(W-3) = \left[\frac{2}{3(W-2)}\right]\left[\frac{1}{3(W-1)}\right]C,
$$

$$
\vdots
$$

$$
h(W-k) - k(W-k-1) = \left[\frac{1 \cdots k}{3^k(W-1)\cdots(W-k)}\right]C.
$$

Simplifying, we have

(22) $\quad h(W-k) - h(W-k-1) = \left[\frac{k!(W-k-1)!}{3^k(W-1)!}\right]C = \frac{C}{3^k\binom{W-1}{k}}.$

To solve for $C$ now, we use a telescoping series to conclude

$$
1 = h(W) - h(0) = C\left[\sum_{k=0}^{W-1}\frac{1}{3^k\binom{W-1}{k}}\right].
$$

TABLE 6

| $x$ | $W = 6$ | $W = 8$ |
|---|---|---|
| 1 | 0.003782 | 0.0004334 |
| 2 | 0.006051 | 0.0006190 |
| 3 | 0.009455 | 0.0008047 |
| 4 | 0.01966 | 0.001139 |
| 5 | 0.08093 | 0.002141 |
| 6 |  | 0.007156 |
| 7 |  | 0.05228 |

Using the last two formulas, we can compute $h(x)$ for each $0 < x < W$, given $W$. Here and in what follows numerical results were obtained by using $C$ and/or Matlab programs. Our Table 6 gives the values of $h(x)$.

Since $P_W(T_W^+ < T_0) = P_{W-1}(T_W < T_0) = h(W-1)$, this gives the values of $a$ quoted in Section 2.2.

(20) implies that $h(X_n)$ is a martingale so if $a < x < b$, then

$$(23) \qquad P_x(T_a < T_b) = \frac{h(b) - h(x)}{h(b) - h(a)}, \qquad P_x(T_a > T_b) = \frac{h(x) - h(a)}{h(b) - h(a)}.$$

From this we can compute Green's function $G_b(x, y) = E_x(N_y(T_b)) =$ the expected number of visits to $y$ starting from $x$ before hitting $b$. If $x, y < b$, then

$$(24) \qquad E_x(N_y(T_b)) = \frac{P_x(T_y < T_b)}{P_y(T_y^+ > T_b)} = \frac{P_x(T_y < T_b)}{p(y, y+1)P_{y+1}(T_y > T_b)}.$$

In words, the numerator gives the probability that we get to $y$ before hitting $b$. If we reach $y$, then we will return a geometric number of times with "success" probability $P_y(T_y^+ > T_b)$. In order to hit $b$ before returning to $y$, we must go to $y + 1$ on the first jump and then hit $b$ before $y$. Some concrete examples that we will need are

$$E_0(N_6(T_7)) = \frac{1}{p(6, 7)} = 12,$$

$$(25)$$

$$E_0(N_5(T_7)) = \frac{1}{p(5, 6)P_6(T_5 > T_7)} = \frac{p(6, 5) + p(6, 7)}{p(5, 6)p(6, 7)} = 80.$$

Our next goal is to compute expected hitting times. We could do this by summing the Green's function: $E_x T_b = \sum_{y < b} G_b(x, y)$. However, we will also need results for conditioned chains, so we will do this by solving equations. If $r(x, y)$ is any irreducible nearest neighbor transition probability on $\{0, 1, \ldots, W\}$, then $u(x) = E_x T_b$ satisfies

$$u(x) = 1 + r(x, x+1)u(x+1) + r(x, x)u(x) + r(x, x-1)u(x-1)$$

Table 7

| x | y = 8 | 7 | 6 | 5 | 4 | 3 | 2 | 1 |
|---|---|---|---|---|---|---|---|---|
| 0 | 69104.23 | 3569.23 | 449.66 | 104.37 | 36.91 | 16.43 | 7.71 | 3.00 |
| 1 | 69101.23 | 3566.23 | 446.66 | 101.37 | 33.91 | 13.43 | 4.71 | 0 |
| 2 | 69096.51 | 3561.51 | 441.94 | 96.66 | 29.20 | 8.71 | 0 | |
| 3 | 69087.80 | 3552.80 | 433.23 | 87.94 | 20.49 | 0 | | |
| 4 | 69067.31 | 3532.31 | 412.74 | 67.46 | 0 | | | |
| 5 | 68999.86 | 3464.86 | 345.29 | 0 | | | | |
| 6 | 68654.57 | 3119.57 | 0 | | | | | |
| 7 | 65535.00 | 0 | | | | | | |

for $0 < x < b$, with $u(b) = 0$ and $u(0) = 1 + r(0,0)u(0) + r(0,1)u(1)$.

Since transition probabilities must sum to 1,

$$r(x, x+1)[u(x) - u(x+1)] = 1 + r(x, x-1)[u(x-1) - u(x)]$$

and it follows that

(26)      $$[u(x) - u(x+1)] = \frac{1 + r(x, x-1)[u(x-1) - u(x)]}{r(x, x+1)}.$$

Using $u(0) - u(1) = 1/r(0,1)$, we can iterate to find the successive differences and then use $u(b) = 0$ to determine the function. This procedure is enough for numerical computations so we will not give a formula for $u$. The values of the hitting times $E_x T_y$ when $W = 8$ are given in Table 7. Note that $u(x) - u(x+1)$ does not depend on $y$, so the differences between values in two successive rows are constant, and are equal to the value sitting above the 0.

We compare $E_0 T_W$ with $E_\pi T_W$ in Section 2.1, where $E_\pi T_W = \sum_x \pi(x) E_x T_W$.

Next we need to compute hitting times for the chain conditioned to hit $W$ before 0. To compute the transition probability,

$$q(x, y) = P(X_1 = y | X_0 = x, T_W < T_0)$$
$$= \frac{P_x(X_1 = y) P_y(T_W < T_0)}{P_x(T_W < T_0)} = \frac{p(x, y)h(y)}{h(x)}.$$

One can use (26) for this chain as well. As Table 8 for $W = 8$ shows, this dramatically reduces the hitting times from the unconditioned values $E_x T_8 \geq 65{,}535$.

To compute $E_0 S$ (recall that $S$ is defined in Approximation 3 in Section 2.3), we begin by noting that $E_0 S = E_0 \tau_{W-2} + E_{W-2} S$. Then by considering what happens on the first jump from $W - 2$,

$$E_{W-2} S = (1 - \rho)[p(W-2, W-1) \cdot 1 + p(W-2, W-2) \cdot E_{W-2} S$$
$$+ p(W-2, W-3) \cdot (E_{W-3} \tau_{W-2} + E_{W-2} S)].$$

Solving, we have

$$E_{W-2}S = (1-\rho)\frac{p(W-2, W-1) + p(W-2, W-3)E_{W-3}\tau_{W-2}}{1-(1-\rho)(1-p(W-2, W-1))}.$$

TABLE 8

| $x$ | $E_x(T_8\|T_8 < T_0)$ | $E_x(T_0\|T_8 < T_0)$ |
|---|---|---|
| 8 | 0 | 47.229002 |
| 7 | 2.156068 | 46.229002 |
| 6 | 8.152877 | 45.138092 |
| 5 | 19.963106 | 43.696338 |
| 4 | 31.794219 | 41.811331 |
| 3 | 39.459771 | 39.184505 |
| 2 | 43.829002 | 35.059746 |
| 1 | 46.229002 | 26.937262 |
| 0 | 47.229002 | 0 |

**Acknowledgments.** We would like to thank Eric Siggia for pointing out the problem in [20]. Andrew D. White visiting professor David Aldous provided helpful advice on the use of the Poisson clumping heuristic. Laurent Saloff-Coste served on Deena Schmidt's Ph.D. committee, tutored us on $\tau_2$ and other aspects of convergence rates of Markov chains. We would like to thank the referee and the associate editor for many detailed comments that helped us improve our paper.

DEPARTMENT OF MATH
523 MALOTT HALL
CORNELL UNIVERSITY
ITHACA, NEW YORK 14853
USA
E-MAIL: rtd1@cornell.edu

CENTER FOR APPLIED MATH
RHODES HALL
CORNELL UNIVERSITY
ITHACA, NEW YORK 14853
USA
E-MAIL: deena@cam.cornell.edu